\documentclass[12pt]{article}

\usepackage{amsgen,amsmath,amstext,amsbsy,amsopn,amsfonts,amssymb}
\usepackage[dvips]{graphicx}

\title{On a new approach to classification of associative algebras}
\author{Vladimir Dergachev}

\def\tr{\textrm{\normalfont tr\,}}

\def\rank{\textrm{\normalfont rank\,}}
\def\rk{\rank}

\def\Mat{\textrm{\normalfont Mat}}

\def\Stab{\textrm{\normalfont Stab}}

\def\Nil{\textrm{\normalfont Nil}}

\def\ker{\textrm{\normalfont ker}\,}
\def\AF{\left.A\right|_F}
\def\LTAF{\left.\tilde{A}\right|_F}
\def\RTAF{\left.\tilde{A}^T\right|_F}

\def\imply{\Rightarrow}
\def\implies{\imply}
\def\intersect{\cap}

\newtheorem{theorem}{Theorem}

\newtheorem{lemma}[theorem]{Lemma}

\newenvironment{proof}{\begin{trivlist}\item[\hskip%
\labelsep{\bf Proof\quad}]}%
{\hfill\qed\end{trivlist}}
\newcommand{\qed}{{\unskip\nobreak\hfil\penalty50\hskip .001pt \hbox{}
          \nobreak\hfil
          \vrule height 1.2ex width 1.1ex depth -.1ex
           \parfillskip=0pt\finalhyphendemerits=0\medbreak}\rm}

\makeatletter
\newcount\@refno
\def\nextref#1#2#3#4{\advance\@refno\@ne
\if@filesw \immediate\write\@auxout
   {\string\bibcite{#1}{\number\@refno}}\fi   
    {{\number\@refno}.\quad{ #2}, {#3}, { #4}.\hfill\\}}

\newcommand{\references}{
\section*{References}
\frenchspacing
\entries\par}

\newcommand{\entries}{
\nextref{DK}
{Vladimir Dergachev and Alexandre Kirillov}
{\em Index of Lie algebras of Seaweed type}
{to appear in {\em Journal of Lie theory}}
\nextref{D}
{Vladimir Dergachev}
{\em Some properties of index of Lie algebras}
{math.RT/0001042}
\nextref
{D1}
{J.Dixmier}
{``Enveloping algebras''}
{American Mathematical Society, Providence (1996) 1-379}
\nextref{E1}{Elashvili, A. G}
{\em Frobenius Lie algebras}
{Funktsional. Anal. i Prilozhen. {\bf 16} (1982), 94--95}
\nextref{E2}{---}
{\em On the index of a parabolic subalgebra 
of a semisimple Lie algebra}
{Preprint, 1990}
\nextref{GK}
{I.M.Gelfand and A.A.Kirillov}
{\em Sur les corps li\'es aux alg\'ebres enveloppantes des alg\'ebres de Lie}
{ Publications math\'ematiques {\bf 31} (1966) 5-20}
\nextref{P}
{R.S. Pierce}
{Associative algebras}
{Springer-Verlag, New York, 1982}

}

\begin{document}
\maketitle
\tableofcontents

\section{Introduction}

The study of associative algebras has usually centered on studying subalgebras,
ideals, homomorphisms and identities - i.e. properties describable within the 
category of algebras and their modules. One of the most famous results include 
the study of commutative algebras by use of multiplicative functionals and 
primitive ideals and the "classification" theorem stating that any finite dimensional
associative algebra is a semi-product of its Jacobson radical and direct
sum of matrix algebras over a skew-field.

Unfortunately, these results do not help in a situaton when the researcher has
a multiplication table for the algebra and wishes to learn as much as possible 
about it.

In this paper we introduce a method by using which one can decompose a finite
dimensional associative algebra into subspaces which provide valuable information
about the structure of the algebra. These subspaces generally will not be
ideals, or even algebras, and their value is in simplifying the operation of
multiplication when considered in basis subordinate to this decomposition.

\section{The method}

We start with a finite dimensional associative algebra $\mathfrak A$, which
we will assume contains unity $1$. Pick any basis $\left\{e_i\right\}$ in $\mathfrak A$.
The matrix $A=\left\{e_ie_j\right\}$ is the multiplication table of $\mathfrak A$
written in the basis $\left\{e_i\right\}$.

Pick a linear functional $F\in \mathfrak A^*$. This choice will affect the 
properties of the resulting decomposition. We define $\left.A\right|_F$
to be the matrix $\left\{F(e_ie_j)\right\}$.

Let 
$$
\ker^L\AF=\left\{x\in\mathfrak A: \forall y \in \mathfrak A \imply A(xy)=0\right\}
$$
$$
\ker^R\AF=\left\{x\in\mathfrak A: \forall y \in \mathfrak A \imply A(yx)=0\right\}
$$
$$
\Nil_F=\ker^L\AF \intersect \ker^R\AF
$$

\begin{theorem}\label{kernels}The spaces $\ker^L\AF$, $\ker^R\AF$ and $\Nil_F$ obey the following
relations with respect to multiplication in $\mathfrak A$:
$$
\begin{array}{clcl}
1. & \ker^L\AF \cdot  \mathfrak A & \subset & \ker^L\AF \\
2. & \mathfrak A \cdot  \ker^R\AF & \subset & \ker^R\AF \\
3. & \ker^L\AF \cdot  \ker^R\AF & \subset & \Nil_F \\
4. & \ker^L\AF \cdot  \Nil_F    & \subset & \Nil_F \\
5. &  \Nil_F    \cdot  \ker^R\AF & \subset & \Nil_F \\
6. & \Nil_F    \cdot  \mathfrak A & \subset & \ker^L\AF \\
7. & \mathfrak A \cdot  \Nil_F & \subset & \ker^R\AF
\end{array}
$$

\end{theorem}
\begin{proof}
1. Let $x\in \ker^L\AF$, $y\in \mathfrak A$. Then for all $z\in \mathfrak A$ we have
$$
F((xy)z)=F(x(yz))=0
$$

2. Let $x\in \mathfrak A$, $y\in \ker^R\AF$. Then for all $z\in \mathfrak A$ we have
$$
F(z(xy))=F((zx)y)=0
$$

3. True because the image should be a subset of both $\ker^L\AF$ and $\ker^R\AF$.

4. A special case of 3.

5. A special case of 3.

6. A special case of 1.

7. A special case of 1.

\end{proof}

Define $N=\dim \mathfrak A$ and $K=\dim \mathfrak A -\dim \Nil_F$. Let us now
change the basis in such a way that the last $N-K$ vectors span $\Nil_F$. 
Upon rewriting $A$ in the new basis we discover that elements in the last
$N-K$ columns and in the last $N-K$ rows are all zero. Let $\left.\tilde{A}\right|_F$
be the submatrix formed by the first $K$ columns and the first $K$ rows.

Since $\Nil_F$ was defined as two-sided kernel of $\AF$ the matrix 
$\lambda \LTAF+\mu \RTAF$ is non-degenerate and can be considered as a matrix
of a bilinear form on $\mathfrak A/\Nil_F$.

We define the characteristic polynomial $\chi_F(\lambda, \mu)$ of $\mathfrak A$ 
as 
$$
\chi_F\left(\lambda, \mu\right)=\det\left(\lambda \LTAF+\mu \RTAF\right)
$$

The characteristic polynomial is a homogeneous polynomial. 
The points $(\lambda, \mu)$ of the projective line in which it vanishes form
the {\em spectrum of $\mathfrak A$}. For convenience we assume the following
parameterization of $\mathbb P^1$:
$$
(\lambda,\mu)=(1, -\alpha)
$$
We now define $\Stab_F(\alpha)=\ker\left(\LTAF-\alpha \RTAF\right)$, with
$\Stab_F(\infty)=\ker\left(\RTAF\right)$.
The sign in front of alpha was chosen so that $\Stab_F(1)$ would correspond 
to the stabilizer of $F$ in Lie algebra $\mathfrak A^L$. 

Since, in general, $\LTAF$ and $\RTAF$ do not commute we should expect the
need to use Jordan spaces. And indeed, even though $\Stab_F(\alpha)$ are 
transversal, their sum is commonly smaller than $\mathfrak A/\Nil_F$.

Since $\chi_F\left(\lambda, \mu\right)$ is non-zero there exists a value $\alpha_0$
in which it does not vanish. We have
$$
\begin{array}{l}
 \LTAF-\alpha \RTAF =\left(\LTAF-\alpha_0 \RTAF\right)-\left(\alpha-\alpha_0\right)\RTAF=\\
\;\;\;\;\;\;\;\;\;=\left(\LTAF-\alpha_0 \RTAF\right)\left(\alpha - \alpha_0\right)\left(\frac{1}{\alpha-\alpha_0}-\left(\LTAF-\alpha_0 \RTAF\right)^{-1}\RTAF\right)
\end{array}
$$

We define $\tilde{V}^k(\alpha)$ to be the $k$-level Jordan space of the operator 
$$
\left(\LTAF-\alpha_0 \RTAF\right)^{-1}\RTAF
$$
corresponding to the eigenvalue $\Lambda=\frac{1}{\alpha-\alpha_0}$. Recall,
that by definition we have
$$
\tilde{V}^0(\alpha)=\ker \left(\frac{1}{\alpha-\alpha_0}-\left(\LTAF-\alpha_0 \RTAF\right)^{-1}\RTAF\right)
$$
and
$$
\tilde{V}^{k+1}(\alpha)=\left\{x:\exists y\in \tilde{V}^k(\alpha)\imply \left(\frac{1}{\alpha-\alpha_0}-\left(\LTAF-\alpha_0 \RTAF\right)^{-1}\RTAF\right)x=y\right\}
$$
Simple transformations show that
$$
\tilde{V}^0(\alpha)=\Stab_F(\alpha)
$$
and 
$$
\tilde{V}^{k+1}(\alpha)=\left\{x:\exists y\in \tilde{V}^k(\alpha)\imply \left(\LTAF-\alpha\RTAF\right)x=\left(\LTAF-\alpha_0 \RTAF\right)y\right\}
$$
with the special case
$$
\tilde{V}^{k+1}(\infty)=\left\{x:\exists y\in \tilde{V}^k(\infty)\imply \RTAF x=\left(\LTAF-\alpha_0 \RTAF\right)y\right\}
$$
Since $\mathfrak A$ is finite-dimensional the spaces $\tilde{V}^k(\alpha)$ stabilize,
we will use $\tilde{V}(\alpha)$ to denote $\tilde{V}^{\dim \mathfrak A}(\alpha)$.
We have 
$$
\mathfrak A/\Nil_F=\oplus_{\alpha} \tilde{V}(\alpha)
$$

We now define $V(\alpha)$  ($V^k(\alpha)$) to be the preimage of $\tilde(V)(\alpha)$
(correspondingly $\tilde{V}^k(\alpha)$) under the projection $\mathfrak A \rightarrow \mathfrak A/\Nil_F$.
Restating the conditions defining $\tilde{V}^k(\alpha)$ we obtain:
$$
V^0(\alpha)=\Stab_F(\alpha)
$$
and
$$
V^{k+1}(\alpha)=\left\{x\in\mathfrak A:\exists y\in V^{k}(\alpha)\forall z\in\mathfrak A \imply F(xz)-\alpha F(zx)=F(yz)-\alpha_0 F(zy)\right\}
$$
with the special case
$$
V^{k+1}(\infty)=\left\{x\in\mathfrak A:\exists y\in V^{k}(\infty)\forall z\in\mathfrak A \imply F(zx)=F(yz)-\alpha_0 F(zy)\right\}
$$
Also, from the Jordan decomposition and definition of $V(\alpha)$ we have
$$
\forall \alpha : \Nil_F \subset V(\alpha)
$$
and 
$$
\mathfrak A/\Nil_F = \oplus_{\alpha} V(\alpha)/\Nil_F
$$
\begin{theorem}\label{V_alpha_independence}
The spaces $V(\alpha)$ do not depend on the choice of $\alpha_0$. Furthermore,
when considering $V(\alpha)$ for a particular $\alpha$ we can use any $\alpha_0\neq\alpha$.
Also for $\alpha\neq\infty$ we can define $V^{k+1}(\alpha)$ as
$$
V^{k+1}(\alpha)=\left\{x:\exists y\in \tilde{V}^k(\alpha)\imply \left(\LTAF-\alpha\RTAF\right)x=\RTAF y\right\}
$$
\end{theorem}
\begin{proof}
Indeed, all statements hold for $V^0(\alpha)=\Stab_F(\alpha)$ by definition. We have for a $k\geq 0$
$$
V^{k+1}(\alpha)=\left\{x:\exists y\in \tilde{V}^k(\alpha)\imply \left(\LTAF-\alpha\RTAF\right)x=\left(\LTAF-\alpha_0 \RTAF\right)y\right\}
$$
Starting with 
$$
\left(\LTAF-\alpha_0 \RTAF\right)y
$$
we obtain
$$
\begin{array}{rcl}
\left(\LTAF-\alpha_0 \RTAF\right)y&=&\left(\LTAF-\alpha_1 \RTAF\right)y-(\alpha_0-\alpha_1)\RTAF y=\\
&=&\left(\LTAF-\alpha_1 \RTAF\right)y+\frac{\alpha_0-\alpha_1}{\alpha-\alpha_1}\left(\LTAF-\alpha\RTAF\right)y-\\
&&-\frac{\alpha_0-\alpha_1}{\alpha-\alpha_1}\left(\LTAF-\alpha_1\RTAF\right) y=\\
&=&\left(\LTAF-\alpha_1 \RTAF\right)y+\frac{\alpha_0-\alpha_1}{\alpha-\alpha_1}\left(\LTAF-\alpha_1\RTAF\right)z-\\
&&-\frac{\alpha_0-\alpha_1}{\alpha-\alpha_1}\left(\LTAF-\alpha_1\RTAF\right) y=\\
&=&\left(\LTAF-\alpha_1 \RTAF\right)\left(y+\frac{\alpha_0-\alpha_1}{\alpha-\alpha_1}z-\frac{\alpha_0-\alpha_1}{\alpha-\alpha_1} y\right)
\end{array}
$$
where $z$ is either zero (when $k=0$) or such that 
$$
\left(\LTAF-\alpha\RTAF\right)y=\left(\LTAF-\alpha_1\RTAF\right)z
$$
(from definition of $V^k(\alpha)$).

For the special case $\alpha=\infty$ we have
$$
\begin{array}{rcl}
\left(\LTAF-\alpha_0 \RTAF\right)y&=&\left(\LTAF-\alpha_1 \RTAF\right)y-(\alpha_0-\alpha_1)\RTAF y=\\
&=&\left(\LTAF-\alpha_1 \RTAF\right)y-(\alpha_0-\alpha_1)\left(\LTAF-\alpha_1\RTAF\right)z=\\
&=&\left(\LTAF-\alpha_1 \RTAF\right)\left(y-(\alpha_0-\alpha_1)z\right)
\end{array}
$$

The last statement is proved similarly. Starting with
$$
V^{k+1}(\alpha)=\left\{x:\exists y\in \tilde{V}^k(\alpha)\imply \left(\LTAF-\alpha\RTAF\right)x=\left(\LTAF-\alpha_0 \RTAF\right)y\right\}
$$
we derive
$$
\begin{array}{rcl}
\left(\LTAF-\alpha_0 \RTAF\right)y&=&\left(\LTAF-\alpha \RTAF\right)y+(\alpha-\alpha_0)\RTAF y=\\
&=&\RTAF z+(\alpha-\alpha_0)\RTAF y=\\
&=&\RTAF \left(z+(\alpha-\alpha_0)y\right)
\end{array}
$$
where $z$ is either $0$ (when $k=0$) or such that $\left(\LTAF-\alpha\RTAF\right)y=\RTAF z$.
In the opposite direction, given 
$$
V^{k+1}(\alpha)=\left\{x:\exists y\in \tilde{V}^k(\alpha)\imply \left(\LTAF-\alpha\RTAF\right)x=\RTAF y\right\}
$$
we compute
$$
\begin{array}{rcl}
\RTAF y&=& \frac{1}{\alpha_0-\alpha}\left(\LTAF-\alpha \RTAF\right)y-\frac{1}{\alpha_0-\alpha}\left(\LTAF-\alpha_0\RTAF\right)y= \\
&=&\frac{1}{\alpha_0-\alpha}\left(\LTAF-\alpha_0 \RTAF\right)z-\frac{1}{\alpha_0-\alpha}\left(\LTAF-\alpha_0\RTAF\right)y=\\
&=&\left(\LTAF-\alpha_0\RTAF\right)\frac{z-y}{\alpha_0-\alpha}
\end{array}
$$
where $z$ is either $0$ (when $k=0$) or is such that 
$$
\left(\LTAF-\alpha\RTAF\right)y=\left(\LTAF-\alpha_0\RTAF\right)z
$$

\end{proof}

\newpage
\section{Properties of the spaces $V^k(\alpha)$}
We will now investigate how the spaces $V^k(\alpha)$ interact with respect to
the operation of multiplication.

\begin{theorem}\label{V_mult} For all $\alpha\neq \infty, \beta\neq \infty$ we have 
$$
V^k(\alpha)\cdot V^m(\beta) \subset V^{k+m}(\alpha \beta)
$$
\end{theorem}
\begin{proof}
First, let us establish this property for $V^0(\alpha)=\Stab_F(\alpha)$ and
$V^0(\beta)=\Stab_F(\beta)$.
We have
$$
\Stab_F(\alpha)=\left\{x_1\in\mathfrak A:\forall z\in\mathfrak A\imply F(x_1z)-\alpha F(zx_1)=0\right\}
$$
and
$$
\Stab_F(\beta)=\left\{x_2\in\mathfrak A:\forall z\in\mathfrak A\imply F(x_2z)-\beta F(zx_2)=0\right\}
$$
For all $z\in\mathfrak A$ we have
$$
\begin{array}{l}
F((x_1x_2)z)-\alpha\beta F(z(x_1x_2))=\\
\;\;\;\;\;\;\;\;=F(x_1(x_2z))-\alpha F((x_2z)x_1)+\alpha F(x_2(zx_1))-\alpha\beta F((zx_1)x_2)=0
\end{array}
$$
and thus $(x_1x_2)\in \Stab_F(\alpha\beta)$.

Now let us apply induction. Assume that for all $k+m<T$ the statement is true
and consider the case $k+m=T$. By theorem \ref{V_alpha_independence} we have
$$
V^k(\alpha)=\left\{x_1\in\mathfrak A:\exists y_1\in V^{k-1}(\alpha):\forall z\in\mathfrak A \imply F(x_1z)-\alpha F(zx_1)=F(zy_1)\right\}
$$
and
$$
V^m(\beta)=\left\{x_2\in\mathfrak A:\exists y_2\in V^{m-1}(\beta):\forall z\in\mathfrak A \imply F(x_2z)-\beta F(zx_2)=F(zy_2)\right\}
$$

For any $z\in\mathfrak A$ we have:
$$
\begin{array}{l}
F(x_1x_2z)-\alpha\beta F(zx_1x_2)=\\
\;\;\;\;\;\;\;=F(x_1(x_2z))-\alpha F((x_2z)x_1)+\alpha F(x_2(zx_1))-\alpha\beta F((zx_1)x_2)=\\
\;\;\;\;\;\;\;=F(x_2zy_1)+\alpha F(zx_1y_2)=\\
\;\;\;\;\;\;\;=\beta F(zy_1x_2)+F(zy_1y_2)+\alpha F(zx_1y_2)=\\
\;\;\;\;\;\;\;=F\left(z\left(\beta y_1x_2 +y_1y_2 +\alpha x_1y_2\right)\right)
\end{array}
$$
By induction hypothesis we have $\left(\beta y_1x_2+y_1y_2+\alpha x_1y_2\right)\in V^{T-1}(\alpha\beta)$
and thus $x_1x_2\in V^T(\alpha\beta)$. 
\end{proof}

\begin{theorem} For all $\alpha\neq 0, \beta\neq 0$ we have 
$$
V^k(\alpha)\cdot V^m(\beta) \subset V^{k+m}(\alpha \beta)
$$
\end{theorem}
\begin{proof}
Let us notice that if we create a new algebra $\mathfrak A^*$ with new multiplication
operation
$$
a * b:=b\cdot a
$$
the space $V_{\mathfrak A}(\alpha)$ corresponding to the algebra $\mathfrak A$
is exactly the space $V_{\mathfrak A^*}(1/\alpha)$ corresponding to the algebra
$\mathfrak A^*$. Thus, the statement of this theorem is a consequence of theorem \ref{V_mult}.
\end{proof}

The symmetry $\alpha \rightarrow 1/\alpha$ is also reflected in the dimensions
of spaces $V(\alpha)$ and $\Stab_F(\alpha)$. Because of the fact that 
$$
\det \left(\lambda \LTAF+\mu\RTAF\right)=\det \left(\lambda \LTAF+\mu\RTAF\right)^T=
\det \left(\lambda \RTAF+\mu\LTAF\right)
$$
and theory of Jordan decomposition (stating that multiplicities of roots of
characteristic polynomial are equal to $\dim V(\alpha)/\Nil_F$) we must have
$$
\dim V(\alpha)/\Nil_F=\dim V(1/\alpha)/\Nil_F
$$
and
$$
\dim V(\alpha)=\dim V(1/\alpha)
$$

Similarly for $\Stab_F(\alpha)$ we have
$$
\begin{array}{l}
\dim \Stab_F(\alpha)=\dim \ker \left(\lambda \LTAF+\mu \RTAF\right)=\\
\;\;\;\;\;\;\;\;\;\;\;\;\;\;\;\;\;\;\;\;\;\;\;\;\dim \ker \left(\lambda \LTAF+\mu \RTAF\right)^T=\\
\;\;\;\;\;\;\;\;\;\;\;\;\;\;\;\;\;\;\;\;\;\;\;\;\dim \ker \left(\lambda \RTAF+\mu \LTAF\right)=\dim \Stab_F(1/\alpha)
\end{array}
$$
due to the well known fact $\rk A=\rk A^T$ of linear algebra.

\section{Multiplicative functionals, ideals and $\AF$}
An important classical method of exploring associative algebras is the study 
of ideals and homomorphisms. Multiplicative functionals (homomorphisms into the
base field) have long played an important role in commutative algebras - especially
in commutative $C^*$-algebras. 

It is easy to see that when $F$ is a multiplicative functional $\AF$ is symmetric
and has rank $1$. It turns out that all functionals with $\rk \AF=1$ are multiplicative.

\begin{theorem}Assume that $F$ is such that $\ker^L \AF=\ker^R \AF=\Nil_F$. Then
$\Nil_F$ is an ideal.
\end{theorem}
\begin{proof}
Indeed this is so because of statements 6 and 7 of theorem \ref{kernels}.
\end{proof}

\begin{theorem}Let $F$ be a functional on $\mathfrak A$ such that $F(1)=1$ and
$\rk \AF=1$. Then $F$ is multiplicative.
\end{theorem}
\begin{proof}
Let us count the dimensions. Since $\rk \AF=1$ it must be that $\dim \mathfrak A/\ker^L\AF=
\dim \mathfrak A/\ker^R\AF=1$ and $\dim \mathfrak A/\Nil_F \le 2$. Also, since
$1$ is always an element of $Stab_F(1)$ and $F(\Nil_F)=0$ it must be that
$\dim \Stab_F(1)/\Nil_F\ge 1$. Yet, since the factor spaces $\ker^L\AF/\Nil_F$,
$\ker^R\AF/\Nil_F$ and $\Stab_F(1)/\Nil_F$ are transversal we must have
$$
\dim \ker^L\AF/\Nil_F+\dim \ker^R\AF/\Nil_F +\dim \Stab_F(1)/\Nil_F \le 2
$$
The only way this is possible is when 
$$
\dim \ker^L\AF/\Nil_F=\dim \ker^R\AF/\Nil_F=0
$$
Thus $\Nil_F$ is an ideal and $\ker F=\Nil_F$. Consequently, as $F(1)=1$, $F$ 
must be multiplicative.
\end{proof}

Note that this theorem holds even for an infinite-dimensional algebra, the 
cause being that, due to the condition $\rk \AF=1$, the bulk of $\mathfrak A$
is residing inside $\Nil_F$ and hence all factor spaces have finite dimension. 

\newpage
\section{Jacobson's radical and spaces $V^k(\alpha)$}
The following famous result (see \cite{??}) describes the structure of finite dimensional associative
algebras:
\begin{theorem}Any finite dimensional associative algebra is a semi-direct
sum of its Jacobson radical and direct sum of matrix rings over skew-fields.
\end{theorem}

Since full matrix rings are well-studied the question arises whether the spaces
$V^k(\alpha)$ provide anything beyond this result. For example, the Jacobson's
radical would be contained within $\Nil_F$ for any $F$ produced by the pullback
from the dual space to the factor of $\mathfrak A$ by its Jacobson's radical.
However, it is easy to see that $F$ must vanish completely on any space $V(\alpha)$,
when $\alpha\neq 1$. As Jacobson's radical is an invariant feature of an associative
algebra there are plenty of functionals which do not vanish on it.
\section{Analysis of matrix algebra}

We will now consider the case when $\mathfrak A=\Mat_n(\mathbb C)$. The dual
space $\mathfrak A^*$ can be modeled as $\Mat_n(\mathbb C)$, with 
$$
F(X)=tr(FX)
$$

Consider the case of generic $F$, i.e. when $F$ is invertible and has distinct 
eigenvalues $\left\{\nu_i\right\}$. Observe that conjugation is an algebra isomorphism -
and thus similar $F$ produce equivalent decompositions. (This is not specific
to $\Mat_n(\mathbb C)$, but common to all algebras). Pick a basis that 
diagonalizes $F$. Let $\left\{e_{ij}\right\}$ be the matrix elements in this basis.
We have:
$$
e_{ij}e_{km}=\delta_{jk}e_{im}
$$

The matrix $F(e_{ij}e_{km})$ has the form:
$$
\begin{array}{ccccccccccccccccc}
& \vline &  & e_{ii}   &  &\vline& & e_{ij}^+ & &\vline & & e_{ij}^- & \cr
\hline
& \vline& \nu_1  & & 0& \vline &&&& \vline \cr
e_{ii} & \vline&  & \ddots  & &\vline& & 0 & &\vline & & 0\cr
& \vline& 0&  & \nu_n &\vline& & & & \vline & \cr
\hline
& \vline& & && \vline & & & & \vline &\nu_{j'} & & 0 \cr
e_{ij}^+ & \vline & & 0 & &\vline & & 0 & &\vline &&\ddots	\cr
& \vline& & && \vline & & & & \vline &0 & &\nu_{j''}\cr
\hline
& \vline& & && \vline &\nu_{i'} & &0 & \vline & \cr
e_{ij}^- &\vline & & 0 & &\vline & &\ddots & &\vline & & 0\cr
& \vline& & && \vline &0 & &\nu_{i''} & \vline & \cr
\end{array}
$$
where $e_{ij}^+$ denote those $e_{ij}$ with $j>i$ and $e_{ij}^-$ denote $e_{ij}$
with $j<i$.

It is easy to see that the $\Nil_F=0$. The characteristic polynomial is 
$$
\chi_F(\lambda,\mu)=(-1)^{\frac{n(n-1)}{2}}(\lambda+\mu)^n\prod_i\nu_i\prod_{i\neq j}{(\lambda\nu_i+\mu\nu_j)}
=(-1)^{\frac{n(n-1)}{2}}\prod_{i,j}{(\lambda\nu_i+\mu\nu_j)}
$$
We find that $V(1)$ is spanned by $e_{ii}$, and $V(\nu_i/\nu_j)=e_{ij}\mathbb C$.

\section{Analysis of $\mathcal B(\mathbb H)$}
We will now consider the algebra $\mathcal B(\mathbb H)$ of bounded operators on countable Hilbertian space.
The dual space of continuous functionals on this algebra is dependent
upon the topology defined on $\mathcal B(\mathbb H)$: for normed topology it consists
of all operators of trace class, with the familiar formula $F(X)=\tr(FX)$. For
weak and strong operator topologies (which are both weaker than normed one)
the dual space is restricted to operators of finite rank. Obviously, in the latter
case, the $\Nil_F$ space will have finite codimension - and the rest of analysis
will follow closely the case of $\Mat_n(\mathbb C)$. 

The case of normed topology is also similar to $\Mat_n(\mathbb C)$. Operators of
trace class are compact. Assume $F$ to have point spectrum with simple multiplicities
and $0$ to be the only value outside of point spectrum. The characteristic polynomial
is not well-defined. However the spaces $V(\alpha)$ are. $V(1)$ is the space
of all bounded operators that commute with $F$. The space $V(\alpha=\nu_i/\nu_j)$, where
$\nu_i, \nu_j$ belong to the point spectrum of $F$, has dimension $1$. 

The remarkable property of this situation is that the only functionals available
are degenerate ones - i.e. those for which $\lambda\AF+\mu\AF^T$ is never invertible.

\section{Regular functionals}
Following \cite{D1} we first prove the following lemma:
\begin{lemma}\label{Transversality lemma}
Let $W$ be a vector subspace of $\mathfrak A$. Fix $\lambda$ and $\mu$. The set $R$ of all $F\in \mathfrak A^*$ such
that $W \intersect \ker \left(\lambda \AF+\mu \AF^T\right) \neq 0$ is closed in
$\mathfrak A^*$.
\end{lemma}
\begin{proof}
Let $e_1,\ldots,e_n$ be the basis of $\mathfrak A$ such that the first $p$ vectors
form the basis of $W$. The system of equations
$$
\left(\lambda \AF + \mu \AF^T\right)\sum_{i=1}^p \epsilon_i e_i=0
$$
is equivalent to
$$
\sum_{i=1}^{p}\epsilon_i\left(\lambda F(e_ke_i)+\mu F(e_ie_k)\right)=0
$$
where $k$ ranges from $1$ to $n$. The matrix of the latter system has entries that 
are linear in $F$. The existence of the non-zero solution is equivalent to requirement
that all $p$-minors vanish. This proves that the set $R$ can be defined as a solution to
the system of polynomial equations. Hence it is closed.
\end{proof}

\begin{theorem}
Let $S$ be a subspace of $\mathfrak A^*$. 
Let $F$ be such that for a fixed pair $\lambda_0, \mu_0$ the dimension of
$\ker\left(\lambda_0 \AF+\mu_0\AF^T\right)$ is the smallest among functionals
in small neighbourhood of $F$ inside affine set $F+S$.

Then for any $G\in S$, any $x\in \ker\left(\lambda_0\AF+\mu_0\AF^T\right)$ and 
any $y\in \ker\left(\mu_0\AF+\lambda_0\AF^T\right)$
we have
$$
G\left(\lambda_0xy+\mu_0yx\right)=0
$$
\end{theorem}
\begin{proof}
Pick any $G\in S$. 
Let $F_\epsilon=F+\epsilon G$. Choose any subspace $W\subset\mathfrak A$ of complementary
dimension that is
transversal to $\ker\left(\lambda_0\AF+\mu_0\AF^T\right)$. The set of $\epsilon$
for which $V_\epsilon=\ker\left(\lambda_0\left.A\right|_{F_{\epsilon}}+\mu_0\left.A\right|_{F_\epsilon}^T\right)$
remains transversal to $W$ is open by lemma \ref{Transversality lemma} and thus
contains a small neighbourhood of $0$. Let $\left\{e_i\right\}$ be the basis of $\mathfrak A$
with first $p$ vectors forming the basis of $W$.

Pick $x\in  \ker\left(\lambda_0\AF+\mu_0\AF^T\right)$. As $W$ and $\ker\left(\lambda_0\left.A\right|_{F_{\epsilon}}+\mu_0\left.A\right|_{F_\epsilon}^T\right)$
form a decomposition of $\mathfrak A$ there exists $w\in W$ such that
$$
x-w \in \ker\left(\lambda_0\left.A\right|_{F_{\epsilon}}+\mu_0\left.A\right|_{F_\epsilon}^T\right)
$$
I.e. for all $i$ 
$$
\lambda_0F_{\epsilon}((x-w)e_i)+\mu_0F_{\epsilon}(e_i(x-w))=0
$$
Hence
$$
\lambda_0F_{\epsilon}(we_i)+\mu_0F_{\epsilon}(e_iw)=\lambda_0F_{\epsilon}(xe_i)+\mu_0F_{\epsilon}(e_ix)
$$
The homogeneous part of this system is the linear system that defines $W\intersect \ker\left(\lambda_0\left.A\right|_{F_{\epsilon}}+\mu_0\left.A\right|_{F_\epsilon}^T\right)$.
As this intersection is trivial the rank of it must be $p$. Noticing that the right hand
part can be derived by substituting $w=x$ (i.e. our system is of the type $Qw=Qx$) we conclude
that there is a unique solution $w_{\epsilon}$. As $F_{\epsilon}$ is linear in $\epsilon$
$w_\epsilon$ is a rational function in $\epsilon$.

Pick now any $y\in\ker\left(\lambda_0\AF+\mu_0\AF^T\right)$. We have:
$$
\lambda_0F_{\epsilon}((x-w_{\epsilon})y)+\mu_0F_{\epsilon}(y(x-w_{\epsilon}))=0
$$
Differentiating by $\epsilon$ produces:
$$
\lambda_0G_{\epsilon}((x-w_{\epsilon})y)+\mu_0G_{\epsilon}(y(x-w_{\epsilon}))-
\left(\lambda_0F_{\epsilon}(w_{\epsilon}'y)+\mu_0F_{\epsilon}(yw_{\epsilon}')\right)=0
$$
After setting $\epsilon=0$ we have:
$$
\lambda_0G(xy)+\mu_0G(yx)=0
$$
\end{proof}

{\bf Corollary 1.} Let $F$ be a functional such that $\Stab_F(\alpha)$ has the smallest
dimension in a neighbourhood of $F$. ($\alpha$ is fixed). Then for all $x\in \Stab_F(\alpha)$
and $y\in\Stab_F(1/\alpha)$ we have:
$$
xy-\alpha yx=0
$$

{\bf Corollary 2.} Let $F$ be a functional with the smallest dimension of $\Stab_F(1)$.
Then $\Stab_F(1)$ is a commutative subalgebra of $\mathfrak A$.

{\bf Corollary 3.} Let $F$ be a functional with the smallest dimension of $\Stab_F(0)$.
Then 
$$
\Stab_F(0)\cdot \Stab_F(\infty)=\left\{0\right\}
$$
In this case $\Nil_F$ is a subalgebra of $\mathfrak A$ with trivial (identically $0$)
multiplication law.
\references

\end{document}